\documentclass[a4paper,12pt]{amsart}

\newcommand{\Canda}[0]{\v{C}anda}

\usepackage[english]{babel}
\usepackage{amsthm}

\newcommand{\Set}[1]{\left\{\, #1 \,\right\}}
\newcommand{\Span}[1]{\langle\, #1 \,\rangle}
\newcommand{\Order}[1]{\lvert #1 \rvert}

\DeclareMathOperator{\Sym}{Sym}
\DeclareMathOperator{\Alt}{Alt}

\renewcommand{\phi}[0]{\varphi}
\renewcommand{\theta}[0]{\vartheta}

\newcommand{\I}{\text{$\mathbf{I}$}}

\newcommand{\Eh}[0]{\hat{\mathcal{E}}}
\newcommand{\balpha}[0]{\breve{\alpha}}
\newcommand{\bbeta}[0]{\breve{\beta}}
\newcommand{\bgamma}[0]{\breve{\gamma}}

\newtheorem{dummy}{Dummy}
\numberwithin{dummy}{section}
\numberwithin{equation}{section}

\newtheorem{theorem}[dummy]{Theorem}

\newtheorem{prop}[dummy]{Proposition}

\newcommand{\etls}[0]{\textsc{etls}}

\newcommand{\HG}[0]{\text{$\Set{1} < H <  G$}}
\newcommand{\KG}[0]{\text{$\Set{1} < K <  G$}}
\DeclareMathOperator{\wreath}{\mathrm{wr}}

\begin{document}

\bibliographystyle{amsplain}

\date{2 December 2003 --- Version 3.03%
}
%\title{Some remarks on the algebraic properties of cryptosystem PGM}

\title[The round functions of cryptosystem PGM]%
{The round functions of cryptosystem PGM generate the  symmetric
group}

\author{A.~Caranti}

\address[A.~Caranti]{Dipartimento di Matematica\\
  Universit\`a degli Studi di Trento\\
  via Sommarive 14\\
  I-38123 Trento\\
  Italy} 

\email{andrea.caranti@unitn.it} 

\urladdr{http://www.science.unitn.it/$\sim$caranti/}

\author{F.~Dalla Volta}

\address[F.~Dalla Volta]{Dipartimento di Matematica e Applicazioni\\
  Edificio U7\\
  Universit\`a degli Studi di Milano--Bicocca\\
  via Bicocca degli Arcimboldi 8\\
  I-20126 Milano\\
  Italy}

\email{francesca.dallavolta@unimib.it}

\urladdr{https://www.unimib.it/francesca-dalla-volta}

\begin{abstract}
 S.~S.~Magliveras  et al.\  have  described symmetric  and public  key
 cryptosystems based  on \emph{logarithmic signatures}  (also known as
 \emph{group bases}) for finite  permutation groups.

 In this paper we
 show that if $G$ is a  nontrivial finite group which is not cyclic of
 order  a  prime,  or the  square  of  a  prime,  then the  round  (or
 encryption) functions of these  systems, that are the permutations of
 $G$  induced by  the \emph{exact-transversal  logarithmic signatures}
 (also  known as  \emph{transversal group  bases}), generate  the full
 symmetric group on $G$.

 This answers a question of S.~S.~Magliveras,
 D.R.~Stinson and Tran van Trung.
\end{abstract}


\thanks{First author partially supported by MIUR-Italy via PRIN 2001012275
 ``Graded  Lie  algebras  and   pro-p-groups  of  finite  width,  loop
 algebras,  and derivations''.  Second  author partially  supported by
 MIUR-Italy via PRIN ``Group theory and applications''.}

\maketitle

\thispagestyle{empty}

\section{Introduction}

S.~S.~Magliveras   has   described   in~\cite{A}   a   symmetric   key
cryptosystem, called  PGM (for  Permutation Group Mappings),  which is
based  on  \emph{logarithmic signatures}  (also  known as  \emph{group
bases})  for   finite  permutation  groups.   

In~\cite{NewApproaches},
S.~S.~Magliveras, D.R.~Stinson  and Tran van Trung  have proposed two
public  key cryptosystems  $\mathrm{MST}_{1}$  and $\mathrm{MST}_{2}$,
which  are based on  logarithmic signatures.   An implementation  of a
symmetric  block  cipher  TST   based  on  these  ideas  is  described
in~\cite{Hor, Impl}.

These  cryptosystems  are  based  on  certain  round  (or  encryption)
functions, the PGM transformations,  which are the permutations on the
set  $\Set{1, 2,  \dots,  \Order{G}}$,  where $G$  is  a finite  group,
induced by
\emph{exact-transversal       logarithmic      signatures} on $G$
(these are also known as
\emph{transversal group bases}; see
Section~\ref{sec:prelim}     for     the    relevant     definitions).

In~\cite{AlgebraicProperties},    S.~S.~Magliveras   and   N.~D.~Memon
studied the  algebraic properties  of the group  generated by  the set
$\Eh$  of   PGM  transformations,  in   particular  investigating  its
size. This is  because a small group here  would make the cryptosystem
weak, and indeed questions about  the size of the corresponding groups
have been asked (and answered) for DES \cite{KRS, CW, WDES}, AES
\cite{WAES}, and other cryptosystems.

S.~S.~Magliveras and  N.~D.~Memon have proved
in~\cite{AlgebraicProperties} that the group is as big
as possible, subject to some restrictions.
\begin{theorem}[Magliveras-Memon]\label{thm:mm}
Let $G$ be a finite non-Hamiltonian group.

Suppose the order of $G$ is different from
\begin{equation*}
  q, 1 + q^{2}, 1 + q^{3}, \dfrac{q^{k} - 1}{q - 1},
  2^{k-1} (2^{k} \pm 1),
  11, 12, 15, 22, 23, 24, 176, 276,
\end{equation*}
where $q$ is a prime power and $k$ is a positive integer.

Then the group $\Span{\Eh}$ generated by $\Eh$ is the full symmetric
group $\Sym(\Order{G})$.
\end{theorem}

(Here a group is said to be \emph{Hamiltonian} when all of its
subgroups are normal.)

S.~S.~Magliveras, D.R.~Stinson and
Tran van Trung suggest in~\cite{NewApproaches} that the above Theorem
may in fact hold in more general circumstances. The goal of this short
note is to show that this is indeed the case. We prove

\begin{theorem}\label{thm:c-dv}
Let $G$ be a nontrivial finite group. Suppose $G$ is not cyclic of
order a prime, or the square of a prime.

Then the group $\Span{\Eh}$  generated by $\Eh$
is the full symmetric group $\Sym(\Order{G})$.
\end{theorem}

In a cyclic group of prime order we
have $\Eh = \emptyset$, so the result does not hold. We
deal with the case of cyclic groups of order the
square of a prime in Section~\ref{sec:psquare}.

The key  to our approach  is an analysis  of some PGM
transformations from  the point of  view of imprimitive  group actions
(Section~\ref{sec:blocks}). We are then able to avoid a call to the
classification  of $2$-transitive groups  (which is  where the  list of
exceptions  in Theorem~\ref{thm:mm} comes  in), obtaining an
elementary proof of Theorem~\ref{thm:c-dv} (Section~\ref{sec:proof}).

We are grateful to Andrea Lucchini for a useful reference.

\section{Preliminaries}
\label{sec:prelim}

In   this  section  we   recall  briefly   the  definitions   we  need
from~\cite{AlgebraicProperties,   NewApproaches},  and  set   up  some
notation for the rest of  the paper. Two convenient references for the
theory  of (permutation)  groups we  use are~\cite{Rob,  Cam}.   For a
positive integer $n$, we write
\begin{equation*}
  \I_{n} = \Set{0, 1, \dots, n-1}
\end{equation*}

Let $G$ be a finite group. Let
\begin{equation}\label{eq:gamma}
  \Set{1} = G_{0} < G_{1} < \dots < G_{s-1} < G_{s} = G
\end{equation}
be a chain of  subgroups of $G$, with $s \ge 2$.  (So there is no such
chain if $G$ is trivial, or if it has prime order.)  An
\emph{exact-transversal  logarithmic signature}  (\etls) for  $G$ with
respect to~\eqref{eq:gamma}  is an $s$-tuple  $\alpha = (\alpha_{1},
\alpha_{2},  \dots,  \alpha_{s})$,   where  each  $\alpha_{i}$  is  a
bijection between  $\I_{\Order{G_{i}:G_{i-1}}}$ and a  complete set of
right coset representatives of $G_{i-1}$ in $G_{i}$, for $i = 1, \dots
s$.  (These are  called \emph{transversal  group bases}  in \cite{Hor,
Impl}.) In  this paper we  will only  need the case  when $s =  2$, so
that~\eqref{eq:gamma} becomes a chain
\begin{equation*}%\label{eq:two}
  \Set{1} < H <  G.
\end{equation*}
Writing $\mu = \Order{H}$ and $\lambda =
\Order{G:H}$ (so that $\lambda \mu = \Order{G}$), we have that 
$\alpha_{1} : \I_{\mu} \to H$ is a bijection, and $\alpha_{2}$
is a bijection between $\I_{\lambda}$ and a complete set of right coset
representatives of $H$ in $G$.

Writing $n = \Order{G}$, an \etls\ $\alpha$ with respect to \HG\
establishes a bijection between 
$\I_{n}$ and $G$, given by
\begin{equation}\label{eq:breve}
  \begin{aligned}
     \balpha :\ \I_{n} &\to G\\
     x &\mapsto 
     \alpha_{1}(x_{1}) \cdot \alpha_{2}(x_{2}),
  \end{aligned}
\end{equation}
where $x$ is written uniquely as $x = x_{2} + \lambda  x_{1}$, with
$x_{2} \in \I_{\lambda}$,  and  $x_{1} \in \I_{\mu}$.  

The map $\I_{n} \to \I_{\mu} \times \I_{\lambda}$, that
maps $x$ to the pair $(x_{1}, x_{2})$, is known
as a \emph{knapsack transformation}~\cite[Def.~2.19]{Hor}. A proof in
Section~\ref{sec:proof} would be slightly smoother using the
(equivalent) knapsack transformation in which the roles of $x_{1}$ and
$x_{2}$ are reversed. We prefer to stick to the conventions
of~\cite{NewApproaches, AlgebraicProperties}, though. 

Once an  \etls\ $\alpha$  is fixed  (see the comments  
in~Subsection~\ref{subsec:alpha}), one may consider the set of
permutations of $\I_{n}$ given by
\begin{equation}\label{eq:Eh}
  \Eh 
  = 
  \Eh_{\alpha}
  =
  \Set{ 
  \balpha \circ \bbeta^{-1} : \I_{n} \to \I_{n}
  \mid
  \text{$\beta$ an \etls\ for $G$}}. 
\end{equation}
(Here and  in the following,  we compose maps left-to-right.)  This is
the  set  of  PGM  transformations  mentioned  in  the  statements  of
Theorems~\ref{thm:mm}~and
\ref{thm:c-dv}. 
Note that if $\gamma$ is a further \etls, we have
\begin{equation*}
   (\balpha \circ \bgamma^{-1})^{-1}
   \circ
   (\balpha \circ \bbeta^{-1})
   =
   \bgamma \circ \bbeta^{-1}.
\end{equation*}
It follows,  as in~\cite{NewApproaches},  that the group  generated by
$\Eh$  also contains all permutations
\begin{equation*}
   \bgamma  \circ \bbeta^{-1} : \I_{n} \to \I_{n},
\end{equation*}
where $\beta, \gamma$ are \etls\ for $G$.

We  write  $\Sym(X)$  (resp.\  $\Alt(X)$) for  the  symmetric  (resp.\
alternating) group  on a set  $X$; we write $\Sym(n)  = \Sym(\I_{n})$,
and similarly for $\Alt$. In particular, $\Eh \subseteq \Sym(n)$.

\section{Imprimitivity}
\label{sec:blocks}

In this section we analyze the permutations
$\balpha \circ \bbeta^{-1} \in \Sym(n)$, where $\alpha$ and $\beta$
are \etls\ of a certain form, from the point of view of
\emph{imprimitive group actions}.
Our arguments apply  to the case when $\alpha$ is  a fixed \etls\ with
respect to  \HG, and  $\beta$ is another  \etls\ with respect  to \HG,
obtained  from  $\alpha$ via  certain  transformations,  which we  now
describe.

We consider the partition of $G$ in
the right cosets of 
$H$, and certain transformation on $G$ that move cosets to
cosets. The  first such transformation  is obtained by  reordering the
coset representatives in $\alpha$. That is, given a permutation $\tau
\in \Sym(\lambda)$, we obtain a new \etls\ $\beta$ by setting
$\beta_{1} = \alpha_{1}$, and then $\beta_{2}(x_{2}) = \alpha_{2}(x_{2}
\tau)$, for $x_{2} \in \I_{\lambda}$. From~\eqref{eq:breve} we have $x
\bbeta = \beta_{1}(x_{1}) \cdot \beta_{2}(x_{2}) = \alpha_{1}(x_{1})
\cdot \alpha_{2}(x_{2} \tau)$. In other words $x
\bbeta = (x \breve{\tau}) \balpha$, where $x
\breve{\tau} = (x_{2} + \lambda x_{1}) \breve{\tau} = x_{2} \tau + \lambda
x_{1}$;  that   is,  $\bbeta  =  \breve{\tau}  \circ
\balpha$.  

All these  transformations $\breve{\tau}$ act  \emph{imprimitively} on
$\I_{n}$. We  recall that if  a group $S$  acts transitively on  a set
$X$, then we say that  $S$ acts \emph{imprimitively} on $X$ (or simply
that $S$ is \emph{imprimitive}) if there is a partition $\mathcal{P}$ of $X$,
called  a \emph{block system},  whose elements,  called \emph{blocks},
satisfy the following properties:
\begin{enumerate}
\item $S$ maps an element of $\mathcal{P}$ onto another element of
$\mathcal{P}$ (it follows in particular that all elements of
$\mathcal{P}$ have the same order);
\item the elements of $\mathcal{P}$ are proper subsets of $X$, containing at
least two elements.
\end{enumerate}
One says  that $S$ \emph{respects}  the block system  $\mathcal{P}$ on
$X$.  (See  for instance~\cite[Sect.~1.9]{Cam}~or \cite[Sect.~7.2]{Rob},
for further details.  We regret that we are using the  term
\emph{block} in a sense that is different by that of
\cite{AlgebraicProperties, NewApproaches}, but the terminology we use
is well established in the context of permutation groups.)

The
transitive group $S$ is said
to be \emph{primitive} if it is not imprimitive. It is an easy fact
that a  $2$-transitive
group is primitive (see~\cite[7.2.4]{Rob}~or
\cite[Theorem~1.7]{Cam}).

In our context,
the  $\breve{\tau}$ act on $\I_{n}$, respecting the block  system on
$\I_{n}$ given by the blocks
\begin{equation}\label{eq:blocks}
  B_{x_{2}} = \Set{x_{2} + \lambda x_{1} : x_{1} \in \I_{\mu}},
\end{equation}
for $x_{2} \in \I_{\lambda}$.
Now note that $\balpha \circ \bbeta^{-1} = 
\balpha \circ \balpha^{-1} \circ \breve{\tau}^{-1} =
\breve{\tau}^{-1}$. We can then forget about $\alpha$ and $\beta$, and
consider only  the $\breve{\tau}^{-1}$. We  call these transformations
$\breve{\tau}^{-1}$  (or the  $\breve{\tau}$, which  is the  same) the
\emph{blockwise permutations} of the block system $B_{i}$.

The second type of transformations arise from permutations within a
single coset. Choose a fixed coset
representative $\alpha_{2}(z_{0})$, for some $z_{0} \in \I_{\lambda}$,
and a fixed element $h \in H$, and 
consider the $\beta$ that coincides with $\alpha$, but for
$\beta_{2}(z_{0}) = h \cdot \alpha_{2}(z_{0})$. We have $x \bbeta = x
\balpha$, except when $x = z_{0} + \lambda x_{1}$, when we have
\begin{align*}
  x \bbeta = \beta_{1}(x_{1}) \cdot \beta_{2}(z_{0})
  = \alpha_{1}(x_{1}) \cdot (h \cdot \alpha_{2}(z_{0}))
  = (\alpha_{1}(x_{1}) \cdot h) \cdot  \alpha_{2}(z_{0}).
\end{align*}
Now, given a  group $H$, the group homomorphism  $H \to \Sym(H)$ given
by  $h \mapsto  (k \mapsto  k \cdot  h)$ is  called  the \emph{regular
representation} of  $H$.  We write  $\tau_{h}$ for the  permutation of
$\I_{\mu}$   induced  by   the  image   of  $h$   under   the  regular
representation, via  the bijection $\alpha_{1} : \I_{\mu}  \to H$. In
other words, for $x_{1} \in \I_{\mu}$ we write $\alpha_{1}(x_{1})
\cdot h = \alpha_{1}(x_{1} \tau_{h})$.  In
this setting, we have $\bbeta =
\breve{\tau}_{z_{0}, h} \circ \balpha$,  where $\breve{\tau}_{z_{0},
h} = \balpha \circ \bbeta^{-1}$ is the identity on all blocks, except
that on the block 
$B_{z_{0}}$ it will act as $(z_{0}  + \lambda  x_{1})
\breve{\tau}_{z_{0}, h} = z_{0} + \lambda (x_{1} \tau_{h})$. 
We call  these transformations  the \emph{regular
permutations} of the block $B_{z_{0}}$. Clearly, they also respect the
block system~\eqref{eq:blocks}.

Note that the combination of the two classes of transformations we
have described go under the name of \emph{monomial transformations}
in~\cite{AlgebraicProperties}. Monomial transformations alone would yield
$1$-transitivity; however, we will be proving a stronger statement in
Section~\ref{sec:proof}.

The third class of transformations  occurs when we obtain $\beta$ from
$\alpha$ by permuting the elements of $H$, that is, by
taking $\beta_{2} = \alpha_{2}$, and then $\beta_{1}(x_{1}) = \alpha_{1}(x_{1}
\tau)$, where $\tau \in \Sym(\mu)$.  This yields, proceeding
as above, transformations of the form $x
\breve{\tau} =  (x_{2} + \lambda x_{1}) \breve{\tau}=  x_{2} + \lambda
(x_{1} \tau)$. In other  words, these \emph{diagonal permutations} act
with the same  permutation at the same time on  all the blocks. (Here
we regard elements in different blocks to be the same if they have the
same $x_{1}$ coordinate.) Again,
these transformations respect the block system~\eqref{eq:blocks}.

\section{Proof of Theorem~\ref{thm:c-dv}}
\label{sec:proof}

\subsection{$2$-transitivity}
\label{subsec:2-trans}

We begin with showing that $\Span{\Eh}$ acts $2$-transitively on $G$.

Fix a nontrivial, proper subgroup $H$  of $G$, and consider
the setting of Section~\ref{sec:blocks}. If  $x, x' \in \I_{n}$ are in
different blocks, and $y, y' \in \I_{n}$ are also in different blocks,
there  is a  composition of  blockwise and  regular  permutations that
carries  $x$  onto $y$  and  $x'$  onto $y'$.  In  fact,  first use  a
blockwise  permutation to  carry $x$  within  the block  to which  $y$
belongs, and  $x'$ within  the block to  which $y'$ belongs.  (To avoid
complicating notation  unnecessarily, we keep  the names $x$  and $x'$
for  the images  of $x$  and $x'$  under this  permutation.)  Then use
regular permutations within  the two blocks to carry  $x$ onto $y$ and
$x'$ onto $y'$.

If $x$ and $x'$ are in the same block $B$, and $y$ and $y'$ are in the
same block $C$, first apply a blockwise permutation to carry $B$ onto
$C$, and then use a diagonal permutation (which induces the full
symmetric group on each block) to carry $x$ onto $y$ and $x'$ onto $y'$.

We are left with the case when $x$ and $x'$ are in the same block $B$,
while   $y$  and   $y'$   are  in   different   blocks.  Clearly   the
transformations of  Section~\ref{sec:blocks} are not  enough here,
as a $2$-transitive group is primitive.

Given the  above, however,  it will  be enough to  find an  element of
$\Eh$  that fixes  $x'$, and  moves  $x$ out  of $B$.   By applying  a
blockwise permutation,  and a  diagonal one, we  may assume that  $B =
B_{0}$,  the zeroth  block of  $\alpha$, and  $x' =  0$.  Suppose thus
$\alpha$  is an  \etls\ with  respect  to \HG,  with $\alpha_{1}(0)  =
\alpha_{2}(0) = 1$, so that   the coset $H \alpha_{2}(0)$
is $H$,  and $x' \balpha =  0 \balpha = 1$.  Write $x \balpha  = h \in
H$. We also  consider another nontrivial, proper subgroup  $K$, and an
\etls\ $\beta$ with respect  to \KG, with $\beta_{1}(0) = \beta_{2}(0)
= 1$.  Let $B_{i}'$ be  the blocks relative  to $\beta$. We  will make
more precise choices  of $H$, $K$ and $\beta$  later, according to the
properties of $G$.

We note  first that since $G$ is  nontrivial, and it is  not cyclic of
order a prime or  the square of a prime, it has  at least two distinct
nontrivial, proper subgroups. Moreover,  if all the nontrivial, proper
subgroups have the same order $p$, then $p$ is a prime number, and so
$G$ is a (non-cyclic) elementary abelian $p$-group of order $p^{2}$.

Accordingly, we distinguish two cases.
Suppose first that $G$ has two nontrivial, proper subgroups $H$ and
$K$, with $\Order{H} < \Order{K}$. We have thus $\Order{B_{i}'} =
\Order{K} > \Order{B_{0}} = \Order{H}$ for all $i$. If $h \in K$, so
that $h \bbeta^{-1}  \in B_{0}'$, we may modify  $\beta$ by a diagonal
permutation, so  that $0 \bbeta =  1$ still holds,  but $h \bbeta^{-1}
\notin B_{0}$,  as $\Order{B_{0}} <  \Order{B_{0}'}$. We have  thus $0
(\balpha \circ  \bbeta^{-1}) = 0$  and $x (\balpha  \circ \bbeta^{-1})
\notin B_{0}$, as requested.  If  $h \notin K$, so that $h \bbeta^{-1}
\in B_{i}'$,  for some $i \ne 0$,  we may modify $\beta$  by a regular
permutation  on $B_{i}'$,  so that  $h \bbeta^{-1}  \notin  B_{0}$, as
$\Order{B_{0}} < \Order{B_{i}'}$. Here, too, we have $0 (\balpha \circ
\bbeta^{-1}) = 0$ and $x (\balpha \circ
\bbeta^{-1}) \notin B_{0}$.

If $G$ is elementary abelian, of order $p^{2}$, let $H$ and $K$ be any
two nontrivial, proper subgroups. Here we have $B_{i} = B_{i}'$ for all
$i$. As $H \cap K = 1$, and $h \ne 1$, we have $h \notin K$, so that
$h \bbeta^{-1} \notin B_{0}$. 

We  have  thus proved    $\Span{\Eh}$  to be  $2$-transitive  in  all
cases. 

\subsection{Completion of the proof} 
Suppose first  that the order $n$  of $G$ is  even, and let $H$  be a
subgroup  of $G$  of  order $2$.

With  respect to  \HG, a  nontrivial
regular permutation  on a given  block will be a  transposition. Since
$\Span{\Eh}$ is $2$-transitive,  it follows that $\Span{\Eh}$ contains
all transpositions, and thus $\Span{\Eh} = \Sym(n)$.

Note that in this case, and with this choice of $H$, the permutations
of Section~\ref{sec:blocks} clearly generate the wreath product $\Sym(2)
\wreath \Sym(n/2)$. This is well-known to be a maximal subgroup of
$\Sym(n)$. We have not used this fact here, but our proof of
$2$-transitivity could be read to mean that $\Span{\Eh}$ contains
properly this wreath product.

When $n$ is odd, we begin with showing that $\Span{\Eh}$ contains a
$3$-cycle. 

Start with  a
diagonal permutation  $\sigma$ which is  the transposition $(a  b)$ on
each block (see the observation at the end of Section~\ref{sec:blocks}). Fix any block $B$.
The regular permutation on the block $B$ induced by a suitable element
of $H$ will be a $p$-cycle of the form $\pi = (a b c \dots)$, for some
$c$. Conjugate  $\sigma$ by $\pi$ to get  a permutation $\sigma^{\pi}
= \pi^{-1} \sigma \pi$
which is the transposition $(a b)$ on all blocks, except that on block
$B$ it will be $(a b)^{(a b c
\dots)} = (b c)$.  We obtain that the product $\sigma^{\pi} \cdot
\sigma$ is  the identity  on all blocks,  except on  block $B$,
where it is  the $3$-cycle $(b c) (a b) = (a b c)$. 

We might now appeal to an observation of C.~Jordan (\cite{Jor},
\cite[Section~5.1, Fact~1]{Cam}) to the effect that a primitive
group  that contains  a $3$-cycle  is  either the  alternating or  the
symmetric group.   We then conclude by observing  that $\Eh$ contains
an  odd   permutation,  and  thus  $\Span{\Eh}  =   \Sym(n)$.   This  follows
from~\cite[Theorem~5.4]{AlgebraicProperties}:  a blockwise permutation
that exchanges  just two blocks will  be the product of  an odd number
$p$ of transpositions.

For completeness, however, we give the short argument (for the case of
$2$-transitive  groups)  that  shows  that $\Span{\Eh}$  contains  all
$3$-cycles, and thus  contains $\Alt(n)$.  Let $a, b,  c$ be any three
distinct  elements   of  $\I_{n}$.   Since   $\Span{\Eh}$  contains  a
$3$-cycle,  and it  is  $2$-transitive,  there are  $d,  e \in  \I_{n}
\setminus  \Set{a,  b,   c}$  such  that  $(b  a  d),   (b  c  e)  \in
\Span{\Eh}$. If  $d = e$,  then $(b c  d) (b a d)^{-1}  = (a b  c) \in
\Span{\Eh}$. If  $d \ne  e$, then $(b  a d)^{(b  c e)} =  (c a  d) \in
\Span{\Eh}$, and we argue as in the previous case.

\subsection{A remark on the choice of $\alpha$}
\label{subsec:alpha}

Concerning   the   definition    of   $\Eh   =   \Eh_{\alpha}$   given
in~\eqref{eq:Eh}, we may note   that all the transformations
$\balpha \circ  \bbeta^{-1}$, that we have considered  in this Section
and the  previous one,  can be  taken with respect  to a  fixed \etls\
$\alpha$,  chosen once  and for  all with  respect to  \HG,  where the
choice of $H$ depends on the  properties of the group, as we have just
seen, and we take (just for simplicity) $\alpha_{1}(0) =
\alpha_{2}(0) = 1$.

\section{The case of the cyclic group of order $p^{2}$} 
\label{sec:psquare}

If $G = \Span{a}$  is a cyclic group of order $p^{2}$,  where $p$ is a
prime,     the     second      part     of     the     argument     of
Subsection~\ref{subsec:2-trans}  does not  work, as  $G$ has  a unique
nontrivial, proper subgroup $H  = \Span{a^{p}}$.  In fact, since $\Eh$
consists      of      the      (imprimitive)      permutations      of
Section~\ref{sec:blocks}, $\Span{\Eh}$ is not $2$-transitive here.

In this case we have the following
\begin{prop}\label{prop:c-dv}
Let $G$ be a cyclic group of order $p^{2}$,  where $p$ is a
prime. Then:
\begin{enumerate}
\item
the group $\Span{\Eh}$
is a proper, imprimitive subgroup of $\Sym(p^{2})$;
\item given  an \etls\ $\alpha$, there exists  a logarithmic signature
  $\gamma$    such   that
\begin{equation*}
\Span{\Eh   \cup    \Set{\balpha   \circ  \bgamma^{-1}}}
=
\Sym(p^{2}).
\end{equation*}
\end{enumerate}
\end{prop}
(The \etls\ $\alpha$ is taken with respect to the only possible choice
\HG.) We refer to  ~\cite{NewApproaches} for the general definition of
\emph{logarithmic signatures}. (These are called \emph{group bases} in
\cite{Hor, Impl}.)   In the  special case of  the cyclic group  $G$ of
order $p^{2}$  we are considering  here, a logarithmic signature  is a
pair of (injective)  maps $\gamma_{1}, \gamma_{2} : \I_{p}  \to G$, so
that   each   element   of   $G$   can  be   written   (uniquely)   as
$\gamma_{1}(x_{1}) \cdot
\gamma_{2}(x_{2})$, for $x_{1}, x_{2} \in \I_{p}$. As in the case of
an \etls, we 
obtain a bijection $\bgamma : \I_{p^{2}} \to G$ as $(x_{2} + p x_{1})
\bgamma = \gamma_{1}(x_{1}) \cdot
\gamma_{2}(x_{2})$.

We choose
$\alpha_{1}(x_{1}) = a^{p x_{1}}$ and $\alpha_{1}(x_{2}) = a^{x_{2}}$,
so that $x \balpha = a^{x}$. Here $B_{0}$ is the set of multiples of
$p$ in $\I_{p^{2}}$. Then we take the logarithmic signature
$\gamma$  defined by $\gamma_{1}(x_{1}) = a^{x_{1}}$ and
$\gamma_{1}(x_{2}) = a^{p x_{2}}$. One sees that $0 (\balpha
\circ \bgamma^{-1}) = 0$, and $p (\balpha
\circ \bgamma^{-1}) = 1$, so that $\balpha
\circ \bgamma^{-1}$ fixes $0$, and takes $p \in B_{0}$ to an element
$1 \notin B_{0}$, as requested.

This  yields  that  the   group  $\Span{\Eh  \cup  \Set{\balpha  \circ
    \bgamma^{-1}}}$ is  $2$-transitive; the rest of  the proof follows
as in Section~\ref{sec:proof}.

%\bibliography{Refs}

\begin{thebibliography}{10}

\bibitem{Cam}
Peter~J. Cameron, \emph{Permutation groups}, London Mathematical Society
  Student Texts, vol.~45, Cambridge University Press, Cambridge, 1999.
  \MR{2001c:20008}

\bibitem{CW}
K.~W. Campbell and M.~J. Wiener, \emph{{DES} is not a group}, Advances in
  Cryptology - Crypto '92 (Santa Barbara, 1992) (Heidelberg) (E.F. Brickell,
  ed.), Lecture Notes in Computer Science, vol. 740, Springer, 1993,
  pp.~512--520.

\bibitem{Impl}
Val{\'e}r \Canda, Tran van Trung, Spyros Magliveras, and Tam{\'a}s Horv{\'a}th,
  \emph{Symmetric block ciphers based on group bases}, Selected areas in
  cryptography (Waterloo, ON, 2000) (Doug Stinson and Stafford Tavares, eds.),
  Lecture Notes in Comput. Sci., vol. 2012, Springer, Berlin, 2001,
  pp.~89--105. \MR{1 895 584}

\bibitem{Hor}
Tam\'as Horv\'ath, \emph{Das {TST}-kryptosystem}, Ph.D. thesis, Fachbereich
  Maschinentechnik der Universit\"at GH Essen, December 1998.

\bibitem{Jor}
C.~Jordan, \emph{Th\'eor\`emes sur les groupes primitifs}, J. Math. Pures Appl.
  (1871), 383--408.

\bibitem{KRS}
Burton~S. Kaliski, Jr., Ronald~L. Rivest, and Alan~T. Sherman, \emph{Is the
  data encryption standard a group? ({R}esults of cycling experiments on
  {DES})}, J. Cryptology \textbf{1} (1988), no.~1, 3--36. \MR{89f:94017}

\bibitem{A}
S.~S. Magliveras, \emph{A cryptosystem from logarithmic signatures of finite
  groups}, Proceedings of the 29th Midwest Symposium on Circuits and Systems
  (Mohammed Ismail, ed.), Elsevier Science Ltd, 1986, pp.~972--975.

\bibitem{NewApproaches}
S.~S. Magliveras, D.~R. Stinson, and Tran van Trung, \emph{New approaches to
  designing public key cryptosystems using one-way functions and trapdoors in
  finite groups}, J. Cryptology \textbf{15} (2002), no.~4, 285--297. \MR{1 944
  653}

\bibitem{AlgebraicProperties}
Spyros~S. Magliveras and Nasir~D. Memon, \emph{Algebraic properties of
  cryptosystem {PGM}}, J. Cryptology \textbf{5} (1992), no.~3, 167--183.
  \MR{93h:94017}

\bibitem{Rob}
Derek J.~S. Robinson, \emph{A course in the theory of groups}, second ed.,
  Graduate Texts in Mathematics, vol.~80, Springer-Verlag, New York, 1996.
  \MR{96f:20001}

\bibitem{WDES}
Ralph Wernsdorf, \emph{The one-round functions of the {DES} generate the
  alternating group}, Advances in cryptology---EUROCRYPT '92 (Balatonf\"ured,
  1992), Lecture Notes in Comput. Sci., vol. 658, Springer, Berlin, 1993,
  pp.~99--112. \MR{94g:94031}

\bibitem{WAES}
Ralph Wernsdorf, \emph{The round functions of {RIJNDAEL} generate the
  alternating group}, Proceedings of the 9th International Workshop on Fast
  Software Encryption, Lecture Notes in Computer Science, vol. 2365,
  Springer-Verlag, Heidelberg, 2002, FSE2002, Leuven, Belgium, February 2002,
  pp.~143--148.

\end{thebibliography}
\providecommand{\bysame}{\leavevmode\hbox to3em{\hrulefill}\thinspace}
\providecommand{\MR}{\relax\ifhmode\unskip\space\fi MR }
% \MRhref is called by the amsart/book/proc definition of \MR.
\providecommand{\MRhref}[2]{%
  \href{http://www.ams.org/mathscinet-getitem?mr=#1}{#2}
}
\providecommand{\href}[2]{#2}

\end{document}